\documentclass[journal]{IEEEtran}
\usepackage{amsfonts}
\usepackage{bm}
\usepackage{amsmath,amsfonts,amsthm,amssymb}
\usepackage{setspace}
\usepackage{Tabbing}
\usepackage{fancyhdr}
\usepackage{lastpage}
\usepackage{extramarks}
\usepackage{chngpage}
\usepackage{algorithmic}
\usepackage{soul,color}
\usepackage{epsfig}
\usepackage{graphicx,float,wrapfig,subfigure}
\usepackage{dsfont}
\usepackage{longtable}
\usepackage{wrapfig}
\usepackage{stfloats}
\usepackage{cite}
\usepackage{graphicx}
\usepackage{array}
\usepackage{multirow}{\tiny }
\usepackage{multicol}
\usepackage{colortbl}
\usepackage{tabularx}
\usepackage{mdwmath}
\usepackage{mdwtab}
\usepackage{color}
\usepackage{verbatim}
\usepackage{amsmath}
\usepackage{tikz}
\usetikzlibrary{calc}
\usepackage{flowchart}
\usetikzlibrary{shapes.geometric, arrows}
\usepackage{overpic}
\usepackage{url}
\usepackage[colorlinks,linkcolor=black,anchorcolor=black,citecolor=black,urlcolor=black]{hyperref} 
\usepackage{breakurl}
\usepackage{footnote}

\allowdisplaybreaks
\begin{document}

\title{
	Joint Planning of PEV Fast-Charging Network and Distributed PV Generation Using the Accelerated Generalized Benders Decomposition
	}

\author{
Hongcai~Zhang,~\IEEEmembership{Student Member,~IEEE,}
Scott~J.~Moura,~\IEEEmembership{Member,~IEEE,}
Zechun~Hu,~\IEEEmembership{Member,~IEEE,}
Wei~Qi,
and~Yonghua~Song,~\IEEEmembership{Fellow,~IEEE}

\thanks{
	\vspace{-3.0mm}
	
	H. Zhang, Z. Hu and Y. Song are with the Department of Electrical Engineering, Tsinghua University, Beijing, 100084, P.~R.~China (email:  zechhu@tsinghua.edu.cn).
	
	S. J. Moura is with the Department of Civil and Environmental Engineering, University of California, Berkeley, California, 94720, USA.

	W. Qi is with the Energy Analysis \& Environmental Impacts Division, Lawrence Berkeley National Laboratory, Berkeley, California, 94720, USA.
	}
	\vspace{-2mm}
}

\maketitle

\begin{abstract} 
	Integration of plug-in electric vehicles (PEVs) with distributed renewable resources will decrease PEVs' well-to-wheels greenhouse gas emissions, promote renewable power adoption and defer power system investments. 
	This paper proposes a multidisciplinary approach to jointly planning PEV fast-charging stations and distributed photovoltaic (PV) power plants on coupled transportation and power networks. 
	First, we develop models of 1) PEV fast-charging stations; 2) highway transportation networks under PEV driving range constraints; 3) PV power plants with reactive power control. 
	Then, we formulate a two-stage stochastic mixed integer second order cone program (MISOCP) to determine the sites and sizes of 1) PEV fast-charging stations; 2) PV power plants. 
	To address the uncertainty of future scenarios, a significant number of future typical load, traffic flow and PV generation curves are adopted. This makes the problem large scale. We design a Generalized Benders Decomposition Algorithm to efficiently solve it.
	To the authors' knowledge, this work is the first that jointly plans both PEV fast-charging stations and PV plants with consideration for PEV driving range limits and reactive PV power control.
	We conduct numerical experiments to illustrate the effectiveness of the proposed method, and validate the benefits of the joint planning and adopting advanced PV reactive power control.
\end{abstract}

\begin{IEEEkeywords}
Plug-in electric vehicle, charging station, PV generation, planning, transportation, AC power flow, second order cone, Accelerated Generalized Benders Decomposition.
\end{IEEEkeywords}

\vspace{-2.0mm}
\section{Introduction}
\IEEEPARstart{I}{ntegration} of PEVs with distributed renewable resources can help reduce PEVs' well-to-wheel greenhouse gas emissions, promote renewable power adoption, alleviate power congestions and defer power system investment.

Encouraging PEVs to consume low-emission renewable power  is one of the key approaches to decarbonizing our modern transportation systems. The emissions of PEVs depend on their energy supply mix. PEVs in areas with high penetration of coal-fired plants may emit more than traditional electric-gasoline hybrid vehicles or even internal combustion engine vehicles\cite{EV_Emission_Argonne2010}. Integrating PEVs with renewable power resources, e.g., wind and PV power etc., can help fully explore PEVs' emission reduction potential whilst promoting renewable power adoption. 

Building PEV charging infrastructure along with distributed renewable power generation to promote local power supplies will also alleviate power congestions, and thereafter, defer power system investments.
The rapidly growing PEV charging power may threaten secure operation of power distribution networks. For \textit{destination charging}, coordinated controlling or vehicle-to-grid technologies can be utilized to alleviate PEV charging power's negative effect \cite{Zhang_V2GTPS2016}, while uncontrollable \textit{fast-charging} power may cause significant power congestions\cite{Zhang_PlanFRLMTSG2016}.
Considering that upgrading distribution systems is usually expensive, installing cheap distributed renewable generation to satisfy congested PEV load is a promising solution.

The growing PEV population is leading to massive investments in charging infrastructure recently. For example, in China, 4.8 million distributed charging spots and more than twelve thousand fast-charging stations are planned for construction by 2020\cite{EVStationPlan_chinadaily}. This investment boom gives the society an opportunity to integrate PEVs with renewable resources at the planning stage, i.e., jointly plan PEV charging stations with distributed renewable resources, so that we can reap the aforementioned benefits.

Among different types of renewable resources, distributed PV power is one of the most promising to supply PEV charging locally because that: 1) They are geographically distributed and close to PEV charging demands; 2) The distribution of PV generation couples with daytime PEV charging power, e.g., workplace charging or fast-charging; 3) Distributed PV generation with advanced grid-connected inverter can help support reactive power control to enhance power quality that may be deteriorated by large-scale PEV integration.

Integrating renewable power with PEV charging stations has been a research hotspot over recent years. 
Most of the published papers focus on economic benefit evaluation or coordinated control strategies. 
Takagi et al. \cite{PEV_PV_Takagi2013} adopted PEV battery-swapping stations to accommodate PV power. 
MacHiels et al. \cite{PEV_PV_MacHiels2014} studied the economic benefit of integrating PV generation with fast-charging stations. 
Brenna et al. \cite{PEV_PV_Brenna2014} and Liao et al. \cite{PEV_PV_Liao2015} demonstrated that coordinated PEV charging could significantly improve distributed PV power integration. 
Alam et al. \cite{PEV_PV_Alam2015} showed that coordinated PEV charging could alleviate voltage rise problems caused by PV power injection.

Few published papers have studied the joint planning of PEV charging stations and renewable power generation. 
Liu et al. \cite{PEV_PV_Plan_Liu2014} studied joint planning of on-site PV generation and battery-swapping stations. The capacities of PV panels, PEV batteries, and number of PEV chargers are optimized at the same time.
Shaaban et al. \cite{PEV_PV_Plan_Shaaban2014} proposed a multi-year multi-objective planning algorithm for uncoordinated PEV parking lots and renewable generation. 
Moradi et al. \cite{Plan_EE_Moradi2015} developed a multi-objective model to optimize the sites and sizes of charging stations and distributed renewable generation. 
Chandra Mouli et al.  \cite{PEV_PV_Plan_ChandraMouli2016} designed a workplace PEV charging station powered by PV generation with vehicle-to-grid technology. 
Quoc et al. \cite{PEV_PV_Plan_Quoc2016} studied the sizing of a PEV charging station powered by commercial grid-integrated PV systems considering reactive power support. The PEV charging station was connected to an inverter that was controlled in three quadrants.
Amini et al. \cite{PEV_PV_Plan_Amini2017} proposed a two-stage approach to simultaneously allocating PEV charging stations with distributed renewable resources in distribution systems.

This paper focuses on joint planning of PEV fast-charging stations and distributed PV power plants. We develop models to determine the sites and sizes of 1) PEV fast-charging stations; 2) PV power plants on coupled transportation and power networks.
The contributions of the proposed method compared to the aforementioned literature are threefold:
\begin{enumerate}
	\item The PEV traffic flows and charging demands are explicitly modeled on transportation networks by the modified capacitated-flow refueling location model (CFRLM) under PEV driving range constraints. By contrast, the aforementioned literature ignored the mobility constraints of PEVs in transportation systems.
	\item This paper considers the new PV power plants with reactive power control technology so that they can help enhance distribution system reliability. By contrast, the aforementioned literature only considered traditional PV power plants. Note that although reference \cite{PEV_PV_Plan_Quoc2016} also considered reactive power control in the planning model, it used the charging station itself rather than the PV power plant to achieve the control. Besides, we use the second order cone programming (SOCP) to describe the power constraints of a PV inverter so that both the active and reactive power can be accurately optimized; by contrast the model of \cite{PEV_PV_Plan_Quoc2016} approximates reactive power based on given active power.
	\item The proposed model is a two stage stochastic MISOCP model, which can be solved by off-the-shelf solvers and the optimality of the solution can be guaranteed. We also design an Accelerated Generalized Benders Decomposition Algorithm to expedite the computation in large scale scenarios. Furthermore, we prove that the algorithm will converge to the optimal solution after a finite number of iterations. By contrast, the aforementioned literature utilized heuristic optimization methods. 
\end{enumerate}
Numerical experiments are conducted to illustrate the effectiveness of the proposed method. The benefits of the joint planning of charging stations with PV power plants and the adoption of PV reactive power control are discussed.

The models of PEV charging stations, transportation networks and PV power generation are formulated in Section II. Section III introduces the MISOCP planning model. The Accelerated Generalized Benders Decomposition Algorithm is given in Section IV. Case studies are described in Section V and Section VI concludes the paper.

\section{Preliminary Models}\label{sec:SRM}
\vspace{-1mm}
\subsection{PEV Charging Station}
\vspace{-1mm}
We adopt the service rate model developed in \cite{Zhang_SOCPEVModel_2016} to describe a PEV charging station's service ability and  model PEV load as a function of the traffic flow visiting a station. 

We assume a set of PEV types, $\mathcal{K}$, with different driving ranges and charging behaviors; PEVs of type $k\in \mathcal{K}$ arrive in a station at location $i$ following a Poisson process with parameter $\lambda_{i,k}$ and requires $T_k$ units of charging time. We let $y_{i,k}^\text{ev}$ denote the number of Poisson arrivals of type $k$ PEVs in $T_k$ units of time in charging station $i$. Therefore, $y_{i,k}^\text{ev}\sim Poisson(T_k\lambda_{i,k}), \forall k\in \mathcal{K}$.
In the station, the PEVs are served on a first-in first-out basis and no arriving PEVs have to wait. \footnote{This assumption is mild and will not significantly affect the quality of the planning results. Interested readers can refer to \cite{Zhang_SOCPEVModel_2016} for detailed discussion.} 
Based on these assumptions, we model a charging station's service ability based on the following service quality criterion:

\textbf{Criterion 1} The probability that any PEV can be charged for at least its required amount of units of time, i.e., $T_k$ for a type $k$ PEV, $k \in \mathcal{K}$, is $\alpha$ or greater. Mathematically, $\text{Pr}(t_{e_{k}}^\text{d}-t_{e_{k}}^\text{a}\geq T_k)\geq \alpha, \forall e_k, \forall k\in \mathcal{K}$, where, $t_{e_{k}}^\text{d}$ is the departure time and $t_{e_{k}}^\text{a}$ is the arrival time of the PEV $e_{k}$.

Criterion 1 is equivalent to the following Criterion 2\cite{Zhang_SOCPEVModel_2016}:

\textbf{Criterion 2} $\text{Pr}( y_i^\text{ev}\leq y_i^\text{cs}) \geq \alpha,\quad y_i^\text{ev}=\sum_{k}{y_{i,k}^\text{ev}},~ y_{i,k}^\text{ev}\sim Poisson(T_k\lambda_{i,k}), \forall k\in \mathcal{K}$, where $y_i^\text{cs}$ is the number of spots.

Each independent Poisson arrival $Poisson(T_k\lambda_{i,k})$ can be approximated by a Normal distribution, i.e., $y_{i,k}^\text{ev}\sim N(T_k\lambda_{i,k}, T_k\lambda_{i,k} )$. Because the sum of different independent Normal random variables is still normally distributed, we have $y_{i}^\text{ev}\sim N(\sum_{k\in K}T_k\lambda_{i,k}, \sum_{k\in K}T_k\lambda_{i,k})$. Then, Criterion 2 is:
\begin{align}
	&\int_{-\infty}^{y_i^\text{cs}}f(y_i^\text{ev})dy_i^\text{ev}=\Phi (\frac{y_i^\text{cs}-\sum_{k\in K}T_k\lambda_{i,k}}{\sqrt{\sum_{k\in K}T_k\lambda_{i,k}}})\geq \alpha,\label{eqn:station5}
\end{align}
where, $f(\cdot)$ and $\Phi(\cdot)$ are respectively the probability and cumulative density function of the normal distribution. 

In practice, the traffic flow passing by one location $i$ may be composed by different types of PEVs from different OD pairs, and only parts of them need charging. Thus, we have: 
\begin{align}
	&\lambda_{i,k}=\sum_{q\in Q_i}{\lambda_{q,i,k} \gamma_{q,i,k}},\label{eqn:station7}
\end{align}
where, $\lambda_{q,i,k}$ is the type $k$ PEV traffic flow on path $q$; $Q_{(i)}$ is the set of paths through node $i$, $q \in Q_{(i)}$. $\gamma_{q,i,k}$ is a binary variable indicating charge choice of type $k$ PEVs on path $q$ at node $i$: $\gamma_{q,i,k}=1$, if they get charged; $\gamma_{q,i,k}=0$, otherwise.

By (\ref{eqn:station5})--(\ref{eqn:station7}), we have the service ability model for a charging station serving $\mathcal{K}$ types of PEVs in an SOCP form:
\begin{align}
	\begin{split}
		y_i^\text{cs}&\geq  \sum_{q\in Q_i}\sum_{k\in K}{T_k\lambda_{q,i,k} \gamma_{q,i,k}}\\
		&+\Phi^{-1}(\alpha) \sqrt{\sum_{q\in Q_i}\sum_{k\in K}{T_k\lambda_{q,i,k} \gamma_{q,i,k}^2} }.\label{eqn:station8}
	\end{split}
\end{align}

The corresponding average PEV charging load is:
\begin{align}
	P_i^\text{ev} = p^\text{sp}\sum_{q\in Q_i}\sum_{k\in K}{T_k\lambda_{q,i,k} \gamma_{q,i,k}},\label{eqn:station9}
\end{align}

in which, $p^{\text{sp}} $ is the rated power of a charging spot.
\vspace{-1mm}

\subsection{Transportation Network}
Driving range limit is the key characteristic of PEVs. Properly modeling this constraint of PEVs on transportation networks enhances the forecasting accuracy of future PEV charging demands. We utilize the modified CFRLM\_SP proposed in our previous work \cite{Zhang_SOCPEVModel_2016} to explicitly consider the driving range constraints of PEVs and adopt time-varying traffic flows to define the locations of the charging stations. 


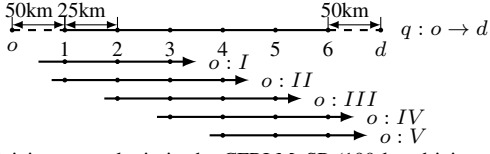
\begin{figure}
	\centering
	\tikzset{>=latex}
	\vspace{-6mm}
	\begin{tikzpicture}[scale=0.7]
	\begin{footnotesize}
	\draw[-,dashed, thick] (-0.,0) -- (1,0)node[right=1.5pt] {};
	\draw[-,dashed, thick] (6,0) -- (7,0)node[right=4.5pt] {$q:o\rightarrow d$};
	\draw[-, thick] (1,0) -- (6,0)node[right=1.5pt] {};
	
	\filldraw [black] (0,0) circle (1pt) node(1)[below=1.5pt]{$o$};
	\filldraw [black] (1,0) circle (1pt) node(2)[below=1.5pt]{1};
	\filldraw [black] (2,0) circle (1pt) node(3)[below=1.5pt]{2};
	\filldraw [black] (3,0) circle (1pt) node(4)[below=1.5pt]{3};
	\filldraw [black] (4,0) circle (1pt) node(5)[below=1.5pt]{4};
	\filldraw [black] (5,0) circle (1pt) node(6)[below=1.5pt]{5};
	\filldraw [black] (6,0) circle (1pt) node(7)[below=1.5pt]{6};
	\filldraw [black] (7,0) circle (1pt) node(8)[below=1.5pt]{$d$};
	
	\draw[->, thick] (1-0.5,-0.6) -- (3.5,-0.6)node[right=1.5pt] {$o:I$};
	\filldraw [black] (1,-0.6) circle (1pt);
	\filldraw [black] (2,-0.6) circle (1pt);
	\filldraw [black] (3,-0.6) circle (1pt);
	
	\draw[->, thick] (0.75,-0.95) -- (4.5,-0.95)node[right=1.5pt] {$o:II$};
	\filldraw [black] (1,-0.95) circle (1pt);
	\filldraw [black] (2,-0.95) circle (1pt);
	\filldraw [black] (3,-0.95) circle (1pt);
	\filldraw [black] (4,-0.95) circle (1pt);
	
	\draw[->, thick] (1.75,-1.3) -- (5.5,-1.3)node[right=1.5pt] {$o:III$};
	\filldraw [black] (2,-1.3) circle (1pt);
	\filldraw [black] (3,-1.3) circle (1pt);
	\filldraw [black] (4,-1.3) circle (1pt);
	\filldraw [black] (5,-1.3) circle (1pt);
	
	\draw[->, thick] (2.75,-1.65) -- (6.5,-1.65)node[right=1.5pt] {$o:IV$};
	\filldraw [black] (3,-1.65) circle (1pt);
	\filldraw [black] (4,-1.65) circle (1pt);
	\filldraw [black] (5,-1.65) circle (1pt);
	\filldraw [black] (6,-1.65) circle (1pt);
	
	\draw[->, thick] (3.75,-2) -- (6.75,-2)node[right=1.5pt] {$o:V$};
	\filldraw [black] (4,-2) circle (1pt);
	\filldraw [black] (5,-2) circle (1pt);
	\filldraw [black] (6,-2) circle (1pt);
	
	\draw[-] (0,0) -- (0,0.25);
	\draw[-] (1,0) -- (1,0.25);
	\draw[<->] (0,0.1) -- (1,0.1) node[above, xshift=-4.75mm] {50km};
	\draw[-] (6,0) -- (6,0.25);
	\draw[-] (7,0) -- (7,0.25);
	\draw[<->] (6,0.1) -- (7,0.1) node[above, xshift=-4.75mm] {50km};
	\draw[-] (1,0) -- (1,0.25);
	\draw[-] (2,0) -- (2,0.25);
	\draw[<->] (1,0.1) -- (2,0.1) node[above, xshift=-4.75mm] {25km};
	
	\end{footnotesize}
	\end{tikzpicture}
	\vspace{-4mm}
	\caption{Driving range logic in the CFRLM\_SP (100 km driving range).}
	\label{fig:sub-path}
	\vspace{-6mm}
\end{figure}

We explain the driving range logic of CFRLM\_SP by Fig.~\ref{fig:sub-path}. A PEV with a driving range of 100 km arrives at node 1 with $D_a=50$ km (it has already traveled 50 km) and needs to depart at node 6 with $D_d=50$ km (so that it can reach its destination after departure). We add pseudo nodes $o$ and $d$ to denote the original node and the destination node respectively and let $d_{o,1}=50$ km and $d_{6,d}=50$ km. The trip setting of our problem is thus equivalent to that a PEV with its battery fully charged leaves at node $o$ and needs to arrive at node $d$ without running out of energy on the road. The travel trajectory of the PEV, i.e., \{$o$, 1, 2, 3, 4, 5, 6, $d$\}, is called a path, i.e., $q$; and a segment of path $q$ is a sub-path. The real nodes on path $q$, i.e., \{1, 2, 3, 4, 5, 6\}, are the candidate locations for charging stations. The \textbf{driving range logic} for a PEV on path $q$ is that any sub-path in $q$ with a distance longer than the PEV's driving range, i.e., 100 km, should cover at least one charging station so that the PEV can travel through path $q$ with adequate charging services.

The driving range logic can be formulated as follows (see Table~\ref{tab:Notation} for additional notation): 
\begin{align}
	&\sum_{i\in \Psi_o^{\text{tn}}}\gamma_{q,i,k} \geq 1, \quad \forall o \in O_{q,k}, \forall q \in Q, \forall k\in K,\label{eqn:1_2}\\
	&\gamma_{q,i,k} \leq x_i^\text{cs}, \qquad \forall q \in Q, i \in \Psi^{\text{tn}}, \forall k\in K,\label{eqn:1_4}\\
	&\underline{y_{i}^\text{cs}}x_{i}^\text{cs} \leq y_i^\text{cs} \leq \overline{y_{i}^\text{cs}}x_{i}^\text{cs}, \qquad \forall i \in \Psi^{\text{tn}}.\label{eqn:1_5}
\end{align}
Equation (\ref{eqn:1_2}) ensures that the PEVs are charged at least once in each sub-path. Equation (\ref{eqn:1_4}) constrains PEVs to charge at nodes with charging stations. Equation (\ref{eqn:1_5}) bounds the number of charging spots.

\vspace{-2mm}
\subsection{PV Generation}
Besides active power generation, PV power plants with fast-reacting and VAR-capable inverters can also generate or consume controllable reactive power which can help enhance reliability of distribution system operations\cite{PV_Q_Turitsyn2011,PV_SDP_Dallranese2014}. The PV generation model with both active and reactive power control can be formulated as an SOCP model as follows:
\begin{align}
&\sqrt{|p^\text{pv}|^2+|q^\text{pv}|^2} \leq \overline{s^\text{pv}},\label{eqn:PV1}\\
&0 \leq p^\text{pv} \leq \overline{p^\text{pv}}, \label{eqn:PV2}\\
&s^\text{pv}=p^\text{pv} + jq^\text{pv},\label{eqn:PV3}
\end{align}
where, $p^\text{pv}$ and $q^\text{pv}$ are respectively the active and reactive power of the PV generation; $\overline{s^\text{pv}}$ is its nameplate apparent power; $\overline{p^\text{pv}}$ is the upper bound of the active power. Equation (\ref{eqn:PV1}) is the constraint for both active and reactive power of the PV generation, which is in the form of an SOCP. The active power is constrained by solar radiation in (\ref{eqn:PV2}). Equation (\ref{eqn:PV3}) calculates the apparent power. In this PV generation model, $q^\text{pv}$ is adjustable and can be either negative or positive.

\vspace{-2mm}
\section{Two-stage Stochastic Joint Planning Model}
We assume the planner is a social planner and has access to the parameters of both the transportation and power systems. It aims to maximize the social welfare. The targeted planning area is a transportation network coupled by a high-voltage distribution system. We also assume that the system can purchase electricity from and sell surplus electricity (at a lower price) to the upper-level power grid.\footnote{Note that this setup is for illustration purposes and not necessarily representative of a particular transportation/power system network. 
}

To represent future probabilistic situations, we first generate a finite set of potential scenarios ($\Omega$), i.e., typical base load, traffic flow and PV generation curves, to represent the future situations for the planning. Then we formulate a two stage stochastic programming model to determine the PEV charging locations, the sites and sizes of both PEV fast-charging stations and PV power plants, and the corresponding distribution system upgrades, i.e., $X=\{\gamma_{q,i,k}, x_i^\text{cs},y_i^\text{cs}, x_m^\text{pv}, \overline{s_{m}^\text{pv}},P_i^\text{sub}\}$.

\vspace{-2mm}
\subsection{Objective}
The planning objective includes the equivalent annual investment costs and the weighted average annual operation costs for all the future scenarios (see Table~\ref{tab:Notation} for the notation):
\begin{align}
Obj&=\min_{{X}}\left\{C^\text{I}\left({X}\right)+\sum_{\omega\in \Omega}{\left(\pi_\omega C^\text{O} \left({X},\omega\right)\right)}\right\}.\label{obj1}
\end{align}

The fist-stage annual investment cost is:
\begin{align}
C^\text{I}\left({X}\right)=&\zeta^\text{cs}\sum_{i\in \Psi^{\text{tn}}}\left(c_{1,i}x_i^\text{cs}+c_{2,i}y_i^\text{cs} +c_{3,i}l_{i}p^\text{sp} y_{i}^\text{cs}+c_{4,i}P_i^\text{sub}\right)\notag\\
&+\zeta^\text{pv}\sum_{m\in \Psi^{\text{dn}}}\left(c_{5,m}x_m^\text{pv}+c_{6,m} \overline{s_{m}^\text{pv}} \right),\label{obj2_1}
\end{align}
where, the substation capacity expansion $P_i^\text{sub}=\max(0,p^\text{sp} y_{i}^\text{cs}-P_{i,0}^\text{sub})$.
The first two terms in the first line of (\ref{obj2_1}) represent the fixed cost of building charging stations and the variable cost in proportion with the number of charging spots. The last two terms in the first line together account for power distribution network upgrade costs, which include the costs for distribution lines and for substation capacity expansion. The two terms in the second line represent the fixed cost per PV plant and the cost per kVA PV panels.

The second stage annual operation costs given the investment decision $X$ for each scenario $\omega$ is:
\begin{align}
	C^\text{O} \left({X},Y_{\omega t},\omega\right)=
	&365\sum_{t} \left( c_\text{e}^{+} p_{0,\omega,t}^{+} \Delta t -c_\text{e}^{-} p_{0,\omega,t}^{-} \Delta t \right)\notag\\
	&+365\sum_{t}\sum_{i\in \Psi^{\text{tn}}} \left( c_\text{p} p_{\text{un},i,\omega,t}^\text{ev} \Delta t \right)\notag\\
	&+365\sum_{t}\sum_{m\in \Psi^{\text{db}}} \sigma\left| v_{m,\omega,t}-v_{0,\omega,t} \right|.\label{obj2_2}
\end{align}
The first two terms in (\ref{obj2_2}) are the system's annual expected energy costs, i.e., the costs for purchasing electricity minus the income by selling surplus electricity. The third term is the penalty for unsatisfied PEV charging power. The fourth term is the penalty for undesirable voltage deviations.\footnote{This term can be easily reformulated as an affine objective by adding two linear inequality constraints for each $\left|\right|$ (absolute value) term.} Coefficient $\sigma$ is used to balance it with the first two monetary objectives.\footnote{In practice, $\sigma$ should be designed according to the system's parameters and the power supply quality requirement. We assume it is given in this paper.}

The second stage optimization variable $Y_{\omega t}$ includes the nodal voltages, branch currents, and PEV charging power etc., which are listed in Table~\ref{tab:Notation}.

\begin{table}
	\renewcommand{\arraystretch}{1.03}
	\vspace{-4mm}
	\begin{footnotesize}
		\caption{Notation Used in the Planning Model}
		\label{tab:Notation}
		\vspace{-1mm}
		\begin{tabular}{p{1.1cm}p{7cm}}
			\hline
			\multicolumn{2}{c}{\textbf{Indices/sets}}\\
			$i/\Psi_{(o)}^{\text{tn}}$&Index/set of transportation nodes (on sub-path $o$), $i\in \Psi_{(o)}^{\text{tn}}$.\\
			$m/n/h$&Index of buses of the distribution network. $m/n/h \in \Psi^{\text{dn}}$. For the substation bus (reference bus), $m/n/h=0$.\\
			$(m,n)/$&Index/set of lines of the distribution network. $(m,n)$ is in  \\
			$\Psi^{\text{db}}$&the order of bus $m$ to bus $n$, i.e., $m\rightarrow n$, and bus $n$ lies between bus $m$ and bus 0. $(m,n) \in \Psi^{\text{db}}$.\\
			$o/O_{(q,k)}$&Index/set of sub-paths (of PEV type $k$ on path $q$), $o \in O_{(q,k)}$.\\
			$\Psi_{(\rightarrow m)}^{\text{dn}}$&Set of buses of the distribution network (that are connected to bus $m$ and bus $m$ lies between them and bus 0).\\
			$\Psi_{m}^{\text{tn}}$& Set of transportation nodes connected to distribution bus $m$.\\
			\hline
			\multicolumn{2}{c}{\textbf{Parameters of the planning model}}\\
			$c_{1,i}$&Fixed costs for building a new station at node $i$, in \$.\\
			$c_{2,i}$&Costs for adding an extra spot in a station at node $i$, in \$.\\
			$c_{3,i}$&Per-unit cost for distribution line at $i$, in \$/(kVA$\cdot$km).\\
			$c_{4,i}$&Per-unit cost for substation capacity expansion at $i$, in \$/kVA.\\
			$c_{5,m}$&Fixed costs for building a PV generation at bus $m$, in \$.\\
			$c_{6,m}$&Costs for adding extra PV panels at bus $m$, in \$/kVA.\\
			$c_\text{e}$&Per-unit cost for energy purchase, in \$/kWh.\\
			$c_\text{p}$&Per-unit penalty costs for unsatisfied PEV power, in \$/kWh.\\
			$\overline{I_{mn}}$&Upper limit of branch current of line $(m,n)$, in kA.\\
			$l_{i}$&Required distribution line length to install a charging station at node $i$, in km.\\
			$N^\text{pv}$&Maximum PV generation number. \\
			$p_{m,\omega,t}^\text{pv,fore}$& Per unit PV power output during $t$ in scenario $\omega$.\\
			$P_{i,0}^\text{sub}$&Initial substation capacity available at node $i$, in kVA.\\
			$s_{m,\omega,t}^\text{b}$&Apparent base load at bus $m$, in kVA.\\
			$\overline{S^\text{pv}}$& Maximum total PV power capacity in the system, in kVA.\\
			$\underline{S_m^\text{pv}}/\overline{S_m^\text{pv}}$& Minimum/maximum PV power capacity at bus $m$, in kVA.\\	$\underline{V_{m}}/\overline{V_{m}}$&Lower/upper limit of nodal voltage at bus $m$, in kV.\\
			$\underline{y_{i}^\text{cs}}/\overline{y_{i}^\text{cs}}$&Minimum/maximum number of charging spots in station $i$.\\
			$Y^\text{ev/pv}$& Service life of the charging stations/PV generation, in year.\\
			$z_{mn}$& Impedance of branch $(m,n)$, in ohm. $z_{mn}^*$ is its conjugate.\\
			$\Delta t$&Time interval, one hour in this paper.\\
			$\zeta^\text{cs/pv}$& Capital recovery factor, which converts the present investment costs into a stream of equal annual payments over the specified time of $Y^\text{cs/pv}$ at the given discount rate $r$. $\zeta=({r(1+r)^{Y^\text{cs/pv}}}) / ({(1+r)_{}^{Y^\text{cs/pv}}-1})$.\\
			$\lambda_{q,i,k,\omega,t}$&Volume of type $k$ PEV traffic flow on path $q$, at node $i$, during time $t$, in scenario $\omega$, in $\text{h}^{-1}$.\\	
			$\pi_\omega$& Probability of scenario $\omega$. \\
			$\omega/\Omega$& Index/set of scenarios. $\omega \in \Omega$.\\		
			\hline
			\multicolumn{2}{c}{\textbf{First stage optimization variables ($X$)}}\\
			$\gamma_{q,i,k}$&Binary charge choice of type $k$ PEVs on path $q$ at node $i$: $\gamma_{q,i,k}=1$, if they get charged; $\gamma_{q,i,k}=0$, otherwise.\\
			$x_i^\text{cs}$&Binary charging station location decision at node $i$: $x_i^\text{cs}=1$, if there is a station at node $i$; $x_i^\text{cs}=0$, otherwise.\\
			$x_m^\text{pv}$&Binary PV generation location decision at bus $m$: $x_m^\text{pv}=1$, if there is PV at bus $m$; $x_m^\text{pv}=0$, otherwise.\\
			$y_i^\text{cs}$&Integer number of charging spots at node $i$.\footnote{To accelerate the optimization speed, we relaxed $y_i^\text{cs}$ to be continuous.}\\
			$P_i^\text{sub}$&Continuous substation capacity expansion at node $i$, in kVA.\\
			$\overline{s_{m}^\text{pv}}$&Continuous invested capacity (maximum nameplate apparent power) of PV panels at bus $m$, in kVA.\\
			\hline
			\multicolumn{2}{c}{\textbf{Second stage optimization variables ($Y_{\omega t}$)}}\\
			$l_{mn,\omega,t}$& Continuous square of the magnitude of line $(m,n)$'s apparent current during $t$ in scenario $\omega$, in $\text{kA}^2$.\\
			$p_{(\text{un},)i,\omega,t}^\text{ev}$&Continuous (unsatisfied) active PEV charging power at node $i$ during $t$ in scenario $\omega$, in kW.\\
			$p_{m,\omega,t}$&Continuous total active power injection at bus $m$ during $t$ in scenario $\omega$, in kW.\\
			$s_{m,\omega,t}$&Continuous total apparent power injection at bus $m$ during $t$ in scenario $\omega$, in kVA. $s_{0,\omega,t}$ (at bus $0$) is also the power consumption of the whole distribution system\cite{OPF_tree_Exactness_Gan2015}.\\
			$s_{m,\omega,t}^\text{ev}$& Continuous apparent PEV power at bus $m$ during $t$ in scenario $\omega$, in kVA.\\
			$S_{mn,\omega,t}$&Continuous apparent power flow from bus $m$ to bus $n$ during $t$ in scenario $\omega$, in kVA.\\
			$v_{m,\omega,t}$&Continuous square of nodal voltage at bus $m$ during $t$ in scenario $\omega$, in kV. Reference voltage $v_{0,\omega,t}$ is given.\\
			$\lambda_{i(,k,\omega,t)}$&Continuous volume of (type $k$) PEVs that require charging at node $i$ (during $t$, in scenario $\omega$), in $\text{h}^{-1}$.\\
			\hline
		\end{tabular}
	\end{footnotesize}
	\vspace{-4mm}
\end{table}

\vspace{-3mm}
\subsection{Constraints}
$\forall i \in \Psi_{}^{\text{tn}}, \forall m \in \Psi^{\text{dn}}, \forall \left(m,n\right) \in \Psi^{\text{db}}, \forall \omega \in \Omega, \forall t :$
\begin{align}
	&\text{service ability constraint of each charging station~}(\ref{eqn:station8}), \\
	&\text{transportation constraints of CFRLM (\ref{eqn:1_2})--(\ref{eqn:1_5})},\\
	&\text{PV power constraints (\ref{eqn:PV1})--(\ref{eqn:PV3})},\\
	&S_{mn,\omega,t}=s_{m,\omega,t}+\sum_{h\in \Psi_{\rightarrow m}^{\text{dn}}}{\left(S_{hm,\omega,t}-z_{hm}l_{hm,\omega,t}\right)},\label{eqn:opf1}\\
	&0=s_{0,\omega,t}+\sum_{h\in \Psi_{\rightarrow 0}^{\text{dn}}}{\left(S_{h0,\omega,t}-z_{h0}l_{h0,\omega,t}\right)}, \label{eqn:opf2}\\
	&v_{m,\omega,t}-v_{n,\omega,t}=2\text{Re}(z_{mn}^{*}S_{mn,\omega,t})- |z_{mn}|^2l_{mn,\omega,t}, \label{eqn:opf3}\\
	&|S_{mn,\omega,t}|^2 \leq l_{mn,\omega,t}v_{m,\omega,t}, \label{eqn:opf4}\\
	&s_{m,\omega,t}=-s_{m,\omega,t}^\text{ev}+s_{m,\omega,t}^\text{pv}-s_{m,\omega,t}^\text{b}, \label{eqn:opf5}\\
	&l_{mn,\omega,t}\leq |\overline{I_{mn}}|^2, \label{eqn:current}\\
	&|\underline{V_{m}}|^2\leq v_{m,\omega,t}\leq |\overline{V_{m}}|^2, \label{eqn:voltage}\\
	&s_{m,\omega,t}^\text{ev}=p_{m,\omega,t}^\text{ev} = \sum_{i\in \Psi_{m}^{\text{tn}}}p_{i,\omega,t}^\text{ev}, \label{eqn:ev1}\\
	&p_{i,\omega,t}^\text{ev} + p_{\text{un},i,\omega,t}^\text{ev}= p^\text{sp}\sum_{q\in Q_i}\sum_{k\in K}{T_k\lambda_{q,i,k,\omega,t} \gamma_{q,i,k}},\label{eqn:ev2}\\
	&\overline{p_{m,\omega,t}^\text{pv}} = {p_{m,\omega,t}^\text{pv,fore}} \overline{s_{m}^\text{pv}}, \label{eqn:PV6}\\
	&\underline{S_m^\text{pv}}x_m^\text{pv}\leq \overline{s_{m}^\text{pv}} \leq \overline{S_m^\text{pv}}x_m^\text{pv}, \label{eqn:PV7}\\
	&\sum_{m\in \Psi^{\text{dn}}} x_m^\text{pv} \leq N^\text{pv}, \label{eqn:PV8}\\
	&\sum_{m\in \Psi^{\text{dn}}} \overline{s_{m}^\text{pv}} \leq \overline{S^\text{pv}}. \label{eqn:PV9}
\end{align}

The branch currents and nodal voltages of the distribution network must satisfy AC power flow constraints (\ref{eqn:opf1})--(\ref{eqn:opf5}) and cannot violate their permitted ranges, i.e., constraints (\ref{eqn:current})--(\ref{eqn:voltage}).\footnote{Though the nodal voltage deviations are already penalized in the objective (\ref{obj2_2}), it is still possible that they may be too large in heavy load scenarios which deteriorates electricity quality significantly. Therefore, we constrain them here to guarantee minimum acceptable electricity quality.} The SOCP relaxation of AC power flow \cite{OPF_tree_Exactness_Gan2015} is adopted. 

We consider hourly power balance in the planning model. The hourly average PEV charging power is calculated by equations (\ref{eqn:ev1})--(\ref{eqn:ev2}).
We assume the base loads $s_{m,\omega,t}^\text{b}$ must be satisfied, while part of the PEV charging power can be discarded due to congestion. When the PEV traffic is low, all PEV charging demands can be fulfilled so that $P_{\text{un},i,\omega,t}^\text{ev}=0$; when the traffic flow grows beyond the system's service ability, some charging demands are not fulfilled and $P_{\text{un},i,\omega,t}^\text{ev}>0$. 

The maximum active PV power constrained in (\ref{eqn:PV6}) depends on the installed PV capacity and the solar irradiation.
Equation (\ref{eqn:PV7}) bounds the capacity of each installed PV power plant. Equations (\ref{eqn:PV8})--(\ref{eqn:PV9}) constrain the total number and the total capacity of the PV power plants in the system, respectively.

The planning model (\ref{obj1})-(\ref{eqn:PV9}) is an MISOCP and can be solved by off-the-shelf solvers, e.g., CPLEX \cite{cplex}.

\vspace{-0mm}
\section{Accelerated Generalized Benders Decomposition Algorithm}
A significant number of scenarios should be considered to enhance planning effectiveness. Thus, the planning model is of high dimension and computationally expensive if directly using off-the-shelf solvers. 
To address this challenge, we adopt the Generalized Benders Decomposition Algorithm\cite{BendersD_Mcdanielt1977}. 

In each scenario, the second stage operation problem solves a 24 hour dynamic optimal power flow problems. However, the corresponding decision variables, e.g., the PEV charging power and the PV generation, in adjacent hours are not coupled. Therefore, when the first stage investment decision, i.e., $X$, is given, the second stage operation problems in every hour of every scenario can be decoupled into low-scale sub-problems that can be efficiently solved in parallel. Based on the above analysis, the proposed algorithm naturally decouples the problem into a master problem, i.e., the planning problem, and a collection of sub-problems, i.e., the operation problem of every hour given $X$. 

For simplicity, we reformulate the original problem (\ref{obj1})-(\ref{eqn:PV9}) into its standard MISOCP form, as follows:
\begin{align}
	\min_{X,Y_{\omega t}} &\quad c^{\top}X+\sum_{\omega\in \Omega}\sum_{t}{d_{\omega t}^{\top} Y_{\omega t}}\label{origin_1}\\
	\text{s.t.:}~& \|A_{\omega tj} X + B_{\omega tj} Y_{\omega t} + e_{\omega tj} \|_2 \leq \notag\\
	&c_{\omega tj}^{\top} X + d_{\omega tj}^{\top} Y_{\omega t} + f_{\omega tj}, \qquad \forall \omega, \forall t, \forall j,\label{origin_2}\\
	& X \in \mathbb{X},\label{origin_3}
\end{align}
where, $wt$ (hour $t$ in scenario $\omega$) is the index of the sub-problems; $j$ is the index of the second order cones; $c$ and $d_{\omega t}$ are objective coefficient vectors; $A_{wtj}$, $B_{wtj}$, $c_{wtj}$, $d_{wtj}$, $e_{wtj}$ and $f_{wtj}$ are respectively coefficient matrices or vectors in the second order cone constraints; $\mathbb{X}$ is the feasible set of $X$ that is irrelevant to sub-problems. Note that parts of $X$ are integer variables, which makes the problem hard to scale.

Given a fixed first stage solution $\hat{X}$, the sub-problem $\omega t$ is a convex SOCP (all the variables are continuous):
\begin{align}
\min_{Y_{\omega t}} &\quad {d_{\omega t}^{\top} Y_{\omega t}}\label{sub_1}\\
\text{s.t.:}~& \|B_{\omega tj} Y_{\omega t} +A_{\omega tj} \hat{X} +  e_{\omega tj} \|_2 \leq \notag\\
&d_{\omega tj}^{\top} Y_{\omega t} + c_{\omega tj}^{\top} \hat{X} + f_{\omega tj}, \quad \forall j. \label{sub_2}
\end{align}

Then, we can obtain the sub-problem's dual problem\cite{OptimizationModel_Ghaoui2014}:
\begin{align}
\max_{\mu_{\omega tj}, u_{\omega tj}, \forall j}& \quad \bigg\{ \sum_{j} -u_{\omega tj}^{\top}\left(A_{\omega tj}\hat{X}+e_{\omega tj}\right) - \notag\\
& \qquad\qquad\qquad \mu_{\omega tj}\left(c_{\omega tj}^{\top} \hat{X} + f_{\omega tj}\right) \bigg\} \label{subdual_1}\\
\text{s.t.:}~& \sum_{j}\left(B_{\omega tj}^{\top} u_{\omega tj} + \mu_{\omega tj} d_{\omega tj}\right) = d_{\omega t}, \label{subdual_2}\\
& \| u_{\omega tj} \|_2 \leq \mu_{\omega tj}, \qquad \forall j, \label{subdual_3}
\end{align}
in which, $\mu_{\omega tj}$ and $u_{\omega tj}$ are the vectors of dual variables. Please refer to Appendix A of the supplementary material for the derivation.

The corresponding master problem is:
\begin{align}
	\min_{X, z} &\qquad c^{\top}X+z \label{master_1}\\
	\text{s.t.:}~&  z \geq \sum_{\omega\in \Omega}\sum_{t} \sum_{j} -\hat{u}_{\omega tj\iota}^{\top}\left(A_{\omega tj}{X}+e_{\omega tj}\right) - \notag\\
	& \qquad \hat{\mu}_{\omega tj\iota}\left(c_{\omega tj}^{\top} {X} + f_{\omega tj}\right), \qquad \iota=1,2,...,\label{master_2}\\
	&X \in \mathbb{X},\label{master_3}
\end{align}
in which, $z$ is an ancillary variable; $\iota$ is the index of iterations. 

The Generalized Benders Decomposition Algorithm solves the master problem (\ref{master_1})-(\ref{master_3}) and the dual of every sub-problem (\ref{subdual_1})--(\ref{subdual_3}) iteratively. In each iteration $\iota$, an optimality cut (\ref{master_2}) is added to the master problem to force its solution to converge to that of the original problem (\ref{origin_1})--(\ref{origin_3}). The algorithm stops when a convergence criterion is met.

We prove that strong duality holds between the sub-problem (\ref{sub_1})--(\ref{sub_2}) and its dual problem (\ref{subdual_1})--(\ref{subdual_3}) in the Appendix B of the supplementary material.
As a result, the cut (\ref{master_2}) in each iteration is always effective before convergence\footnote{If the new cut did not force the master problem to obtain a new solution, then the $LB$ and $UB$ in Table \ref{tab:ABD} are qual so that the solution is optimal.}, and the algorithm will converge to the global optimal solution after a finite number of iterations\cite{BendersD_Mcdanielt1977}.

We utilize two techniques to accelerate the algorithm:
\subsubsection{Relaxing the service ability constraint (\ref{eqn:station8})} Constraint (\ref{eqn:station8}) has no second stage decision variables but should be satisfied for every hour in every scenario (because of different traffic flows). However, it will be binding only in peak traffic hours in practice.\footnote{If the constructed charging spots can satisfy peak-hour traffic flows' charging demands, they can also satisfy the demands during other periods.} Therefore, we relax constraint (\ref{eqn:station8}) as follows:
\begin{align}
\begin{split}
y_i^\text{cs}&\geq  \sum_{q\in Q_i}\sum_{k\in K}{T_k\lambda_{q,i,k,\hat{\omega t}_i} \gamma_{q,i,k}}\\
&+\Phi^{-1}(\alpha) \sqrt{\sum_{q\in Q_i}\sum_{k\in K}{T_k\lambda_{q,i,k,\hat{\omega t}_i} \gamma_{q,i,k}^2} }, ~~ \forall i \in \Psi_{}^{\text{tn}},\label{eqn:station8_R}
\end{split}
\end{align}
where, $\hat{\omega t}_i$ is the index of the sub-problem that has the highest traffic flow at location $i$. We then add constraint (\ref{eqn:station8_R}) directly to the master problem and remove constraint (\ref{eqn:station8}) from every sub-problem. This  approach leads to two benefits: 1) the scale of each sub-problem decreases significantly; 2) the modified sub-problem only solves an optimal power flow problem that allows load shedding which is strictly feasible given any $X$ so that we need not consider feasibility cuts.\footnote{Note that, if constraint (\ref{eqn:station8}) is not relaxed and should be satisfied in every sub-problem, it may be violated given some myopic $X$. As a result, we should add extra iterations to generate feasibility cuts to the master problem.} 

\subsubsection{Relaxing the integer constraints of the master problem}
The master problem is computationally intensive for each iteration, since it contains a significant number of integer variables. We first relax its integer constraints and solve the problem (with higher efficiency) until convergence. Then, we add the integer constraints back to the master problem and conduct extra iterations until the new problem converges. Note that this approach will not affect the optimal solution because the feasible set of the original master problem is a subset of the relaxed master problem. Thus, the optimality cuts generated for the latter is also valid for the former\cite{BendersD_Costa2012}.

The pseudo-code of the algorithm is shown in Table~\ref{tab:ABD}. $\varepsilon_1$ and $\varepsilon_2$ are respectively the relevant gaps at convergence of the original problem and its  relaxed continuous form. 

\begin{table}
	\renewcommand{\arraystretch}{1.05}
	\vspace{-6mm}
	\begin{footnotesize}
		\caption{Accelerated Generalized Benders Decomposition}
		\label{tab:ABD}
		\vspace{-2mm}
		\begin{tabular}{p{0.35cm}p{7.85cm}}
			\hline
			$01$&\textbf{Initialization:} Set iteration number $\iota=0$, lower bound $LB=-\infty$, upper bound $UB=+\infty$, relevant gap $Gap=+\infty$, $flag=0$.\\
			$02$&\textbf{While} termination criteria, i.e., $Gap \leq \varepsilon_2$, not fulfilled, \textbf{do}\\
			$03$&\begin{adjustwidth}{3mm}{0cm}$\iota=\iota+1$.\end{adjustwidth}\\
			$04$&\begin{adjustwidth}{3mm}{0cm}\textbf{Step 0} If $Gap \leq \varepsilon_1$ and $flag=0$, let $UB=+\infty$, $flag=1$.
				\end{adjustwidth}\\
			$05$&\begin{adjustwidth}{3mm}{0cm}\textbf{Step 1} If $flag=1$, solve master problem (\ref{master_1})--(\ref{master_3}); otherwise, solve the relaxed continuous form of (\ref{master_1})--(\ref{master_3}). Update the solution $\hat{X}$ and $\hat{z}$. Let $LB=c^T \hat{X}+\hat{z}$.\end{adjustwidth}\\
			$06$&\begin{adjustwidth}{3mm}{0cm}\textbf{Step 2} Solve each sub-problem's dual problem (\ref{subdual_1})--(\ref{subdual_3}), and update each solution $\hat{u}_{\omega tj\iota}$ and $\hat{\mu}_{\omega tj\iota}$. Let $UB=\min \bigg\{UB,~c^T \hat{X} +\sum_{\omega\in \Omega}\sum_{t}\sum_{j}\bigg(-\hat{u}_{\omega tj\iota}^{\top}\left(A_{\omega tj}\hat{X}+e_{\omega tj}\right) - \hat{\mu}_{\omega tj\iota}$$\left(c_{\omega tj}^{\top} \hat{X} + f_{\omega tj}\right)\bigg)\bigg\}$.\end{adjustwidth}\\
			$07$&\begin{adjustwidth}{3mm}{0cm}\textbf{Step 3} Add a new cut (\ref{master_2}) for iteration $\iota$ to the master problem (\ref{master_1})-(\ref{master_3}).\end{adjustwidth}\\
			$08$&\begin{adjustwidth}{3mm}{0cm}\textbf{Step 4} $Gap=100\% \times (UB-LB)/UB$.\end{adjustwidth}\\
			$09$ &\textbf{End while}\\
			$10$ &\textbf{Output} $\hat{X}$ as the solution.\\
			\hline
		\end{tabular}
	\end{footnotesize}
	\vspace{-4mm}
\end{table}

\section{Case Studies and Conclusions}
\subsection{Case Overview and Parameter Settings}
We consider a 25-node transportation network coupled with a 14-node 110 kV high voltage distribution network to illustrate the proposed planning method. Due to limited space, the detailed parameters of the distribution and the transportation networks are omitted in this paper, but can be found in \cite{Zhang_SOCPEVModel_2016}. 

Seventy-two representative scenarios, i.e., three types of weather (rainy, cloudy, sunny) in weekday and weekend of twelve months, of hourly base load, traffic flow and PV power profiles are generated based on PG\&E load profiles\cite{Loadprofile_PGE}, the National House Travel Survey data\cite{NHTS}, and the National Solar Radiation Data Base\cite{PVprofile_PGE}. 

We assume there are four types of PEVs on the road with equal market share, and their driving ranges per charge are respectively 200, 300, 400 and 500 km. The rated charging power $p^\text{sp}$ is 44 kW, and the average service time to charge the four types of PEVs with empty batteries are about 42, 63, 84, 105 minutes. We also assume $D_a = 100$ km, $D_d=100$ km for all PEVs, $\underline{y_{i}^\text{cs}}=0$, $\overline{y_{i}^\text{cs}}=200$ and $\alpha=80\%$.
The costs of PEV charging stations $c_{1,i}=\$ 163,000$ and $c_{2,i}=\$ 31,640$; the distribution line cost $c_{3,i}=120$ \$/(kVA$\cdot$km). The line distance $l_{i}$ is assumed to be 10\% of the distance between the PEV charging station and its nearest 110 kV distribution node. The substation expansion cost $c_{4,i}=788$ \$/kVA. We assume each original 25 transportation node has 1 MVA surplus substation capacity. 
The electricity purchase cost $c_\text{e}^+=0.094$ \$/kWh \cite{Plan_Zhang2015} and the selling price $c_\text{e}^-$ is 30\% lower. The per-unit penalty cost for unsatisfied charging demand $c_\text{p}=10^{3}$ \$/kWh.
We assume all the nodes (except node 1) in the distribution network are candidate PV locations. The PV generation investment cost $c_{5,m}=0$ \$/VA, $c_{6,m}=1,770$ \$/kVA \cite{PV_Costs_Chung2015}. We also assume that $\sigma= \$10^{-4}$, $Y^\text{cs/pv}=15$, $r=8\%$, $\overline{S^\text{pv}}=90~\text{MVA}$, $\underline{S_m^\text{pv}}=0~\text{MVA}$, $\overline{S_m^\text{pv}}=\infty ~\text{MVA}$, $\forall m$.\footnote{Note that there is usually enough land available in highway networks to build PV power plants. Therefore, we do not limit the $\overline{S_m^\text{pv}}$ here.}

We design six cases, with different PEV traffic flows and maximum numbers of PV power plants with or without reactive power control to illustrate the proposed planning method. The parameters of different cases are illustrated in Table~\ref{tab:case}.

\begin{table}
	\renewcommand{\arraystretch}{1.05}
	\vspace{-6mm}
	\centering
	\begin{footnotesize}
		\caption{The parameters of different cases}
		\vspace{-2mm}
		\begin{tabular}{cccc}
			\hline
			\multirow{2}*{Case}&Max. total number/capacity &Reactive   &Daily PEV\\
			& of PV power plants& power control & traffic flow\\
			\hline
			1&0/0 MVA&--&20000\\
			2&5/90 MVA&No&20000\\
			3&5/90 MVA&Yes&20000\\
			4&0/0 MVA&--&40000\\
			5&5/90 MVA&No&40000\\
			6&5/90 MVA&Yes&40000\\
			\hline 
		\end{tabular}\label{tab:case}
	\end{footnotesize}
	\vspace{-3mm}
\end{table}

We set $\varepsilon_1=0.5\%$, $\varepsilon_2=2\%$ in Table~\ref{tab:ABD} and use CPLEX\cite{cplex} to solve the master problem and sub-problems on a workstation with a 12 core Intel Xeon E5-1650 processor and 64 GB RAM. 

\begin{figure}
	\centering
	\vspace{-1.mm}
	\includegraphics[width=0.6\columnwidth]{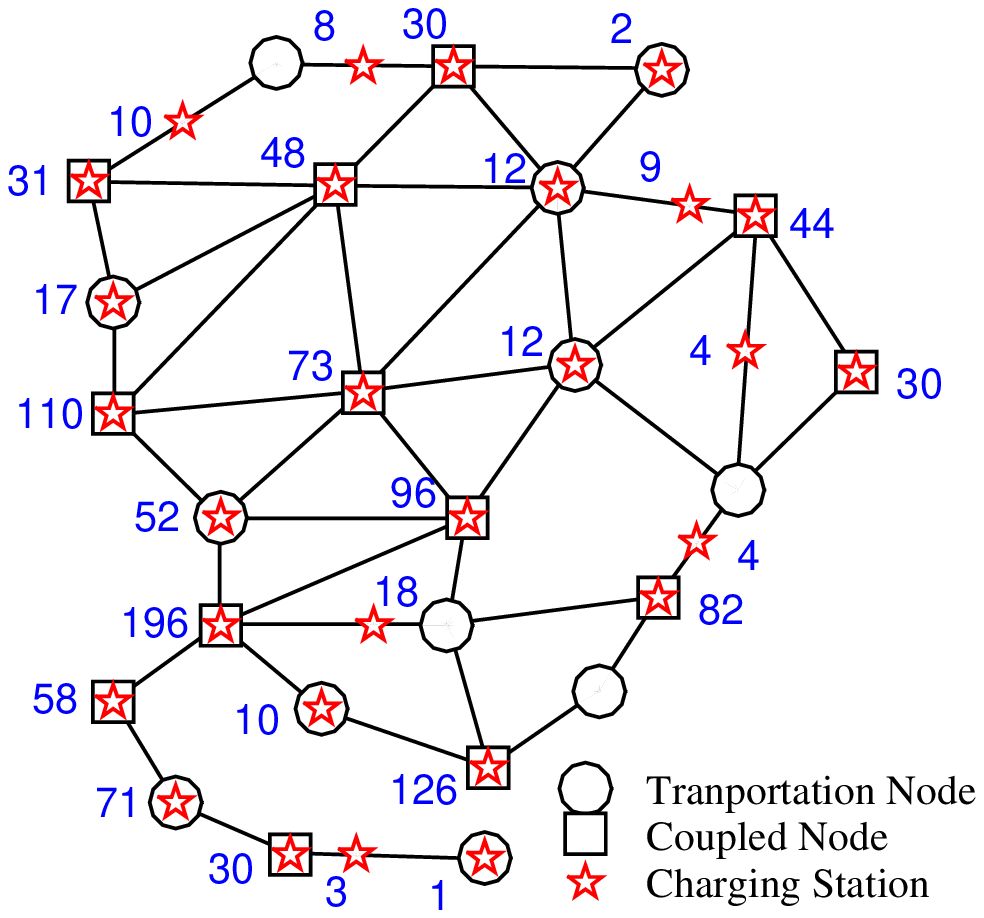}
	\vspace{-7mm}
	\caption{Sites and sizes of PEV charging stations in Case 1. The number next to each station is its capacity, i.e., number of spots.}
	\label{fig:station1}
	\vspace{-5mm}
\end{figure}

\subsection{Planning Results and Analysis}

\begin{table*}[!]
	\renewcommand{\arraystretch}{1.025}
	\vspace{-4mm}
	\centering
	\begin{footnotesize}
		\caption{Summary of the planning results in different cases}
		\vspace{-2mm}
		\begin{tabular}{cccccccccccc}
			\hline
			\multirow{2}*{Case}&Station&Spot &PV &PV capacity&\multicolumn{3}{c}{Investment costs (M\$/year)} & Energy costs & Total costs& Unsatisfied & Solution\\
			& no.& no. &no. & (MVA)& PEV Station & Grid upgrade & PV Plants&  (M\$/year) & (M\$/year) & PEV load (\%) & time (h)\\
			\hline
			1&33&1210&0&0.0&5.39&5.23&0.0&37.94&48.56&0.0&0.5\\
			2&26&1187&4&71.52& 5.17 & 4.58 &14.79&22.47&47.01&0.0&16\\
			3&28&1187&4&73.20& 5.21& 4.58&15.14 &21.95 &46.87&0.0&15\\
			4&44&2279&0&0.0& 9.73& 14.24&0.0 &48.39 &72.37 &1.85 &1.8 \\
			5&31&2287&5&90& 9.50& 11.24&18.61 &28.74&68.50& 0.0 &18\\
			6&30&2285&4&90& 9.48& 11.55&18.61 &28.53 &68.17 &0.0 &18\\
			\hline 
		\end{tabular}\label{tab:result}
	\end{footnotesize}
	\vspace{-2mm}
\end{table*}

\begin{figure*}[!]
	\centering
	\vspace{-1mm}
	\begin{minipage}[t]{2\columnwidth}
		\centering
		\includegraphics[width=0.8\columnwidth]{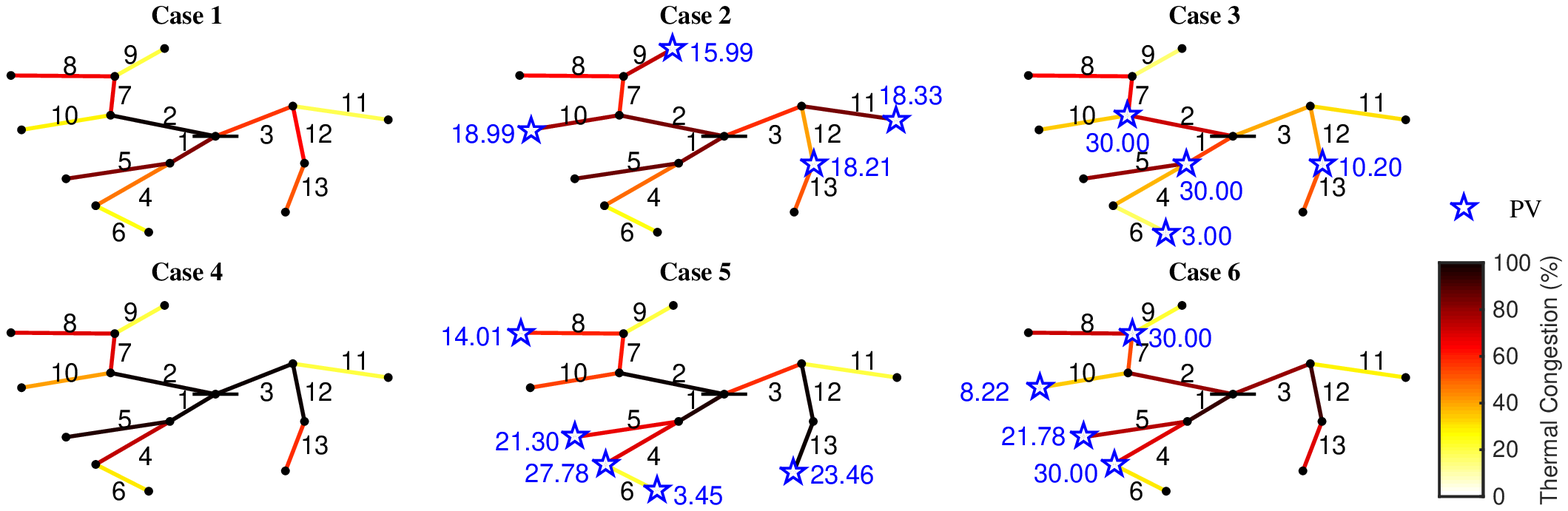}
		\vspace{-4mm}
		\caption{Sites and sizes of PV plants and maximum distribution line thermal congestion levels. The number next to each PV plant is its capacity, in MVA.}
		\label{fig:PV}
	\end{minipage}
	\begin{minipage}[t]{2\columnwidth}
		\centering
		\includegraphics[width=0.88\columnwidth]{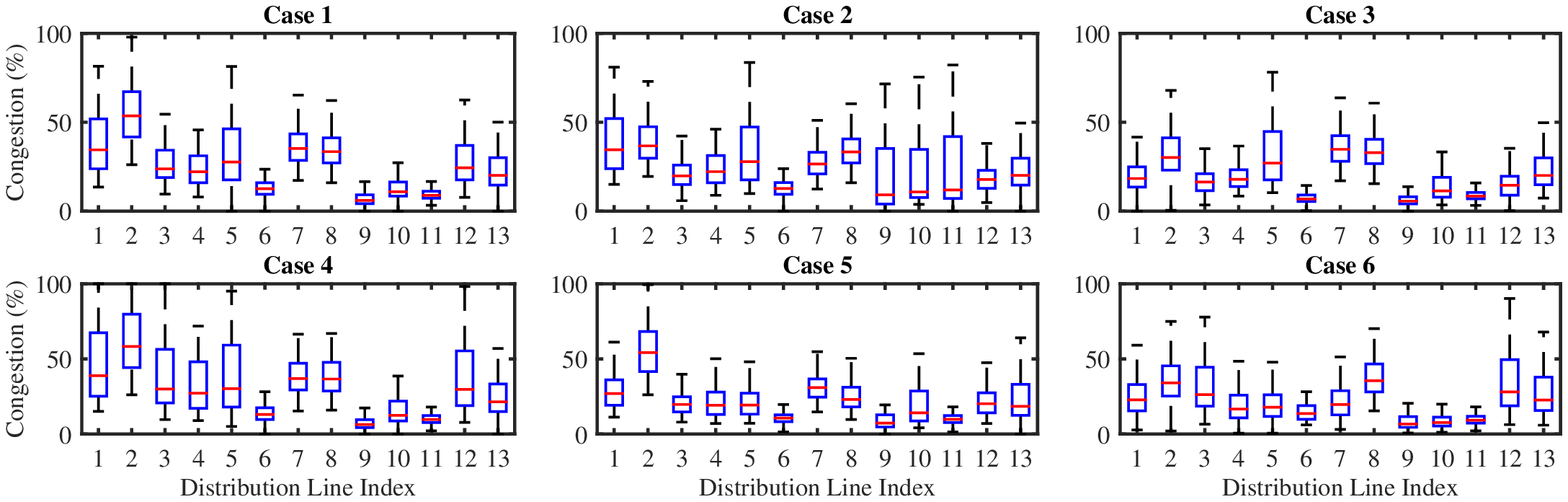}
		\vspace{-3mm}
		\caption{Boxplot of distribution line thermal congestion levels. Line 2 is typically the most congested.}
		\label{fig:congestion}
	\end{minipage}
	\begin{minipage}[t]{2\columnwidth}
		\centering
		\includegraphics[width=0.88\columnwidth]{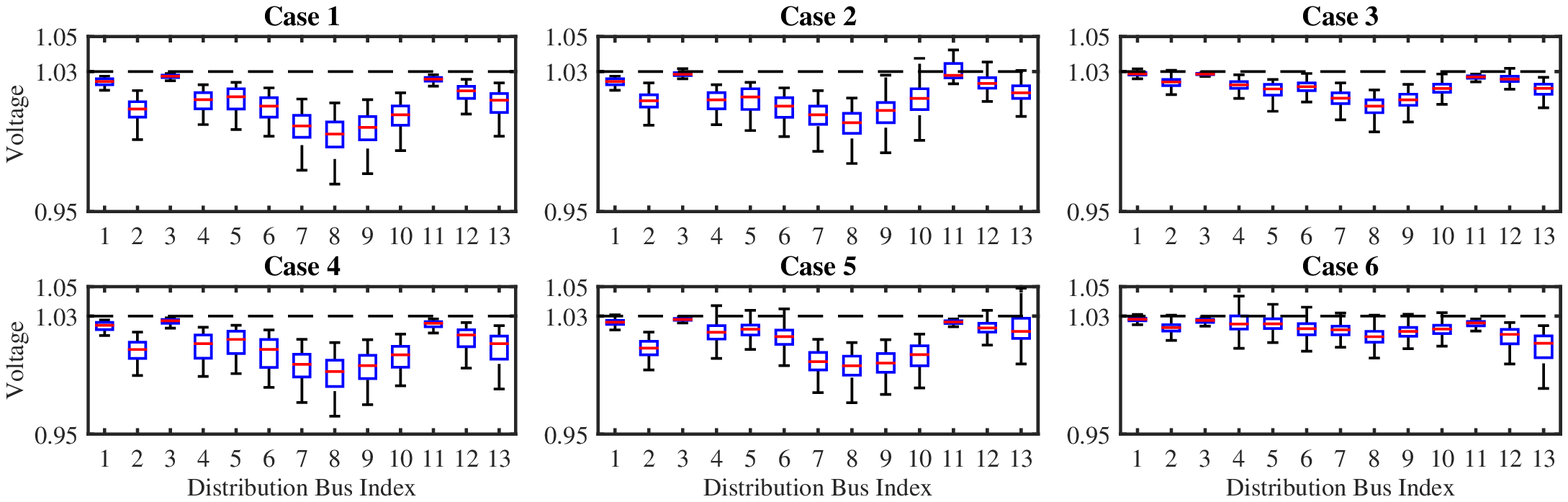}
		\vspace{-3mm}
		\caption{Boxplot of voltages. The reference voltage is 1.03 at reference bus 0. }
		\label{fig:voltage}
	\end{minipage}
	\vspace{-4mm}
\end{figure*}

The summary of the planning results for the six cases are given in Table~\ref{tab:result}. The locations and capacities of PEV charging stations in Case 3 are given in Fig.~\ref{fig:station1} for demonstration. The PV generation and their capacities in different cases are illustrated in Fig.~\ref{fig:PV}. The ratio of a line's current to its thermal capacity, i.e., $100\% \times {\sqrt{l_{mn}}/\overline{I_{mn}}}$, represents its thermal congestion level. The maximum congestion level, i.e., $100\% \times\max_{\omega,t}{\left(\sqrt{l_{mn,\omega,t}}/\overline{I_{mn}}\right)}$, of each distribution line in the six cases are depicted by Colorbars in Fig.~\ref{fig:PV}. The distributions of the line congestion levels and nodal voltages in all the $24\times 72$ hours are respectively illustrated in Figs.~\ref{fig:congestion}--\ref{fig:voltage}.

\subsubsection{Computational efficiency}When jointly planning both PEV charging stations and PV power plants, the scale of the problem is larger; as a result, the solution time is also longer. However, the proposed algorithm can still solve the problems in acceptable time, i.e,. about 18 hours. We can also see that the solution time is longer when the PEV population is larger. That is because larger PEV population leads to higher charging demands and more binding power flow constraints. As a result, the feasible set of the problem is smaller and the algorithm has to conduct more iterations to converge.

\subsubsection{The direct financial benefit for saving total cost}
The planning results show that by jointly building PEV charging stations and PV power plants, the total cost of the system is cut down: the total cost in Case 2 is reduced by 3.19\% compared to Case 1 and the total cost in Case 5 is reduced by 5.35\% compared to Case 4. Though the equivalent annual investment cost is increased, the installed PV power plants generate and sell electricity to the power grid, which significantly decreases the operational costs. 

By utilizing distributed PV generation to supply power locally, the planner has larger flexibility to build PEV charging stations. Compared to Case 1 and Case 4, the overall investment costs on PEV charging stations and the corresponding power grid upgrades in both Case 2 and Case 5 are reduced. This phenomenon is especially prominent in heavy load scenarios. We can observe that in Case 4, much more charging stations are installed than in Case 5. Because some parts of the distribution system are congested, the planner has to build more charging stations elsewhere with higher costs to avoid the PEVs being charged at congested areas. 

The total PV generation capacity and the direct financial  benefit of integrating PEV charging stations with PV generation increase as the PEV population (or load) increases. 

\subsubsection{The indirect benefit by deferring power system investment}
Figs.~\ref{fig:PV}--\ref{fig:congestion} show that investing distributed PV generation can significantly ease distribution line congestion, and therefore, defer power system investment. In Case 2, line 2 is the only one that is congested, which reflects the bottleneck of the system. In Case 4, several distribution lines' capacity constraints are binding, and as a result, 1.87\% of the PEV charging demands cannot be satisfied. By contrast, in the cases with PV generation, no line is congested. Without building new PV power plants, the planner has to upgrade the congested distribution lines (line 2 would be the first choice), which would be much more expensive.

\subsubsection{The benefit of utilizing reactive power control}
By adopting reactive power control for PV generation, the system has larger operational flexibility. As a result, the total cost and the voltage deviations of the system are reduced. Though the monetary benefits seems to be insignificant (less than 1\%'s total cost reduction), Fig.~\ref{fig:voltage} shows that the system with reactive power control has much lower voltage deviations so that it can provide higher quality electricity to customers. Note that, in both Case 2 and Case 4, we can observe significant voltage rises caused by inverse PV power flow. By contrast, in both Case 3 and Case 6, the voltage rises are mild. This advantage will also be much more pronounced at heavy load and high PV penetration scenarios when voltage drops and rises will significantly deteriorate the power quality.


\section{Conclusion}\label{section:conclusion}
We first develop a two-stage stochastic SOCP for jointly planning PEV fast-charging stations and distributed PV power plants on coupled transportation and power networks. Then, we design a Generalized Benders Decomposition Algorithm to efficiently solve the program by decoupling it into a mixed-integer linear master problem and a set of convex SOCP sub-problems.
Our experiments show that investing in distributed PV power plants with PEV charging stations has multiple benefits, e.g., reducing the greenhouse gas emission, promoting renewable power integration, alleviating power congestion caused by large-scale integration of PEVs and thereafter deferring power system investments. The benefits become more prominent when utilizing PV generation with reactive power control, which can help enhance power supply quality.


\bibliographystyle{ieeetr}
\bibliography{ref}

\end{document}



\title{
	Supplementary Material for ``Joint Planning of PEV Fast-Charging Network and Distributed PV Generation Using the Accelerated Generalized Benders Decomposition''
}

\author{
Hongcai~Zhang,~\IEEEmembership{Student Member,~IEEE,}
Scott~J.~Moura,~\IEEEmembership{Member,~IEEE,}
Zechun~Hu,~\IEEEmembership{Member,~IEEE,}
Wei~Qi,
and~Yonghua~Song,~\IEEEmembership{Fellow,~IEEE}

\thanks{
	
	H. Zhang, Z. Hu and Y. Song are with the Department of Electrical Engineering, Tsinghua University, Beijing, 100084, P.~R.~China (email:  zechhu@tsinghua.edu.cn).
	
	S. J. Moura is with the Department of Civil and Environmental Engineering, University of California, Berkeley, California, 94720, USA.
	
	W. Qi is with the Energy Analysis \& Environmental Impacts Division, Lawrence Berkeley National Laboratory, Berkeley, California, 94720, USA.
	}
	\vspace{-4mm}
}

\maketitle


\appendices
\section{The Dual Problem}
\subsection{The Full Formulation of the Sub-problem}
Given a fixed first stage solution $\hat{X}$, the sub-problem $\omega t$ is a convex SOCP (all the variables are continuous), as follows:
\begin{align}
	\min_{Y}~
	&\bigg\{\left( c_\text{e}^{+} p_{0 }^{+} \Delta t -c_\text{e}^{-} p_{0 }^{-} \Delta t \right)+\sum_{i\in \Psi^{\text{tn}}} \left( c_\text{p} p_{\text{un},i }^\text{ev} \Delta t \right) \notag \\
	&+\sum_{m\in \Psi^{\text{db}}} \sigma v_m^\text{d} \bigg\},\label{obj3}\\
\text{s.t.:}~~&\forall i \in \Psi_{}^{\text{tn}}, \forall m \in \Psi^{\text{dn}}, \forall \left(m,n\right) \in \Psi^{\text{db}}: \notag\\
&\sqrt{|p_m^\text{pv}|^2+|q_m^\text{pv}|^2} \leq \overline{s_m^\text{pv}},\label{eqn:3_1}\\
&0 \leq p_m^\text{pv} \leq \overline{p_m^\text{pv}}, \label{eqn:3_2}\\
&s_m^\text{pv}=p_m^\text{pv} + jq_m^\text{pv},\label{eqn:3_3}\\
&S_{mn }=s_{m }+\sum_{h\in \Psi_{\rightarrow m}^{\text{dn}}}{\left(S_{hm }-z_{hm}l_{hm }\right)},\label{eqn:3_4}\\
&0=s_{0 }+\sum_{h\in \Psi_{\rightarrow 0}^{\text{dn}}}{\left(S_{h0 }-z_{h0}l_{h0 }\right)}, \label{eqn:3_5}\\
&v_{m }-v_{n }=2\text{Re}(z_{mn}^{*}S_{mn })- |z_{mn}|^2l_{mn }, \label{eqn:3_6}\\
&|S_{mn }|^2 \leq l_{mn }v_{m }, \label{eqn:3_7}\\
&s_{m}=-s_{m}^\text{ev}+s_{m}^\text{pv}-s_{m}^\text{b}, \label{eqn:3_8}\\
&s_{m}^\text{ev}=p_{m}^\text{ev} = \sum_{i\in \Psi_{m}^{\text{tn}}}p_{i}^\text{ev}, \label{eqn:ev1}\\
&p_{i}^\text{ev} + p_{\text{un},i}^\text{ev}= p^\text{sp}\sum_{q\in Q_i}\sum_{k\in K}{T_k\lambda_{q,i,k} \gamma_{q,i,k}},\label{eqn:ev2}\\
&0 \leq l_{mn }\leq |\overline{I_{mn}}|^2, \label{eqn:3_9}\\
&|\underline{V_{m}}|^2\leq v_{m }\leq |\overline{V_{m}}|^2. \label{eqn:3_10}\\
&v_m^\text{d} \geq v_{m }-v_{0 },\label{eqn:3_11}\\
&v_m^\text{d} \geq -v_{m }+v_{0 },\label{eqn:3_12}\\
&p_{i}^\text{ev}\geq 0, \quad p_{\text{un},i}^\text{ev}\geq 0.\label{eqn:3_13}
\end{align}
in which $v_m^\text{d}$ is the nodal voltage deviation compared to the reference $v_0$.
The objective (\ref{obj3}) is linear; constraints (\ref{eqn:3_1}) and (\ref{eqn:3_7}) are second order cones; the other constraints are all affine. The decision variables is $Y=\{l_{mn }, p_{i}^\text{ev}, p_{\text{un},i}^\text{ev}, s_{0 }, s_{m}^\text{pv}, s_{m}^\text{EV}, S_{mn},  v_{m }, v_m^\text{d}, \forall i \in \Psi_{}^{\text{tn}}, \forall m \in \Psi^{\text{dn}}, \forall \left(m,n\right) \in \Psi^{\text{db}} \}$. We let $\mathcal{D}$ denote the domain of the sub-problem (\ref{obj3})--(\ref{eqn:3_13}), i.e., the intersection of the domains of the objective and the constraint functions of (\ref{obj3})--(\ref{eqn:3_13}). It's obvious that $\mathcal{D}=\mathbf{R}^{d} = \text{relint}~\mathcal{D}$ ($d$ is the dimension of $Y$).

\subsection{The Sub-problem's Dual Problem}
For simplicity, we reformulate the sub-problem (\ref{obj3})--(\ref{eqn:3_13}) in its standard form:
\begin{align}
	p^{*} = &\min_{Y}  ~ {d^{\top} Y}\label{sub_1}\\
	\text{s.t.:}~& \|B_{j} Y +A_{j} \hat{X} +  e_{j} \|_2 \leq d_{j}^{\top} Y + c_{j}^{\top} \hat{X} + f_{j}, \quad \forall j, \label{sub_2}
\end{align}
in which $p^{*}$ is the primal objective.

We follow the procedure in \cite{OptimizationModel_Ghaoui2014} to obtain its dual problem. First, we have
\begin{align}
	p^{*} = &\inf_{Y} \sup_{\mu \geq 0} ~ {d^{\top} Y}  + \sum_j \mu_j \left(\|B_{j} Y +A_{j} \hat{X} +  e_{j} \|_2  \right. \notag\\
	& \qquad\qquad\left. -\left(d_{j}^{\top} Y + c_{j}^{\top} \hat{X} + f_{j}\right)\right)\\
	=&\inf_{Y} \sup_{\|u_j\|_2 \leq \mu_j, \forall j} ~ {d^{\top} Y}  + \sum_j \left(-u_j^{\top} \left(B_{j} Y +A_{j} \hat{X} +  e_{j} \right)  \right. \notag\\
	& \qquad\qquad\qquad~~\left. -\mu_j \left(d_{j}^{\top} Y + c_{j}^{\top} \hat{X} + f_{j}\right)\right),
\end{align}
where we have used the dual representation of the Euclidean norm. $\mu_{j}$ is the dual variable vector of each second order cone and $u_{j}$ is the dual variable vector of each Euclidean norm. 

Then, adopting the max-min inequality\cite{OptimizationModel_Ghaoui2014}, we have 
\begin{align}
d^{*} = &\sup_{\|u_j\|_2 \leq \mu_j, \forall j} \inf_{Y}  ~ {d^{\top} Y}  + \sum_j \left(-u_j^{\top} \left(B_{j} Y +A_{j} \hat{X} +  e_{j} \right)  \right. \notag\\
& \qquad\qquad\qquad~~\left. -\mu_j \left(d_{j}^{\top} Y + c_{j}^{\top} \hat{X} + f_{j}\right)\right), \label{dmax}
\end{align}
which makes $d^{*} \leq q^{*}$.

Solving (\ref{dmax}) for variable $Y$, we obtain the dual problem:
\begin{align}
	d^{*} = & \sup_{\mu_{j}, u_{j}, \forall j}~ \bigg\{ \sum_{j} -u_{j}^{\top}\left(A_{j}\hat{X}+e_{j}\right) - \mu_{j}\left(c_{j}^{\top} \hat{X} + f_{j}\right) \bigg\} \label{subdual_1}\\
	\text{s.t.:}~& \sum_{j}\left(B_{j}^{\top} u_{j} + \mu_{j} d_{j}\right) = d, \label{subdual_2}\\
	& \| u_{j} \|_2 \leq \mu_{j}, \qquad \forall j, \label{subdual_3}
\end{align}
which is still a convex SOCP.

\section{Proof of Strong Duality}
\subsection{The Slater's Condition}
The Slater's Condition provides a sufficient condition for strong duality. We give a brief introduction for it in this part. 

For a convex optimization problem:
\begin{align}
p^{*} = &\min_{x}  ~ f_0(x)\label{convex1}\\
\text{s.t.:}~& f_i (x) \leq 0, \qquad i =1, ..., m,\\
& h_i (x) = 0,\qquad i =1, ..., q,\label{convex3}
\end{align}
we still let $\mathcal{D}$ denote the domain of the problem. Then, we have the following proposition:

\textbf{Proposition 1 (Slater's conditions for convex programs)} \textit{Let $f_i, i=0, ..., m$, be convex functions, and let $h_i, i=0, ..., q$, be affine functions. Suppose further that the first $k \leq m$ of the $f_i$ functions, $i=1,..., k$, are affine (or let $k=0$, if none of the $f_i, i=0, ..., m$, is affine). If there exists a point $x\in \text{relint}~\mathcal{D}$ such that 
\begin{align}
&f_i (x) \leq 0, \qquad i =1, ..., k,\\
&f_i (x) < 0, \qquad i =k+1, ..., m,\label{inequality}\\
& h_i (x) = 0, \qquad i =1, ..., q,
\end{align}
then strong duality holds between the primal problem (\ref{convex1})--(\ref{convex3}) and its dual problem. Moreover, if the primal problem is bounded, i.e., $p^{*}>-\infty$, then the dual optimal value equals to the primal optimal value.}\cite{OptimizationModel_Ghaoui2014}

In the following section, we will use the above proposition to prove strong duality of the sub-problem (\ref{sub_1})--(\ref{sub_2}) and its dual problem (\ref{subdual_1})--(\ref{subdual_3}). We say an inequality constraint to be ``strictly satisfied'' to refer to that it is ``satisfied with strict inequality'' as (\ref{inequality}).

\subsection{Proof of Strong Duality}
We assume that the system can be operated without PV generation and PEV charging power, and the constraints of nodal voltages of the distribution system is not binding. Note that this is a very mild assumption, because the distribution system is usually operated with the voltage deviations being well controlled. Otherwise, the power quality is poor and extra voltage control devices should be installed for the system. 

We first let $s_{m}^\text{pv} =0$ and $s_{m}^\text{EV} =0$, $\forall m \in \Psi^{\text{dn}}$. With constraints (\ref{eqn:ev1})--(\ref{eqn:ev2}), we can directly calculate variables $p_{i}^\text{ev}=0$ and $p_{\text{un},i}^\text{ev}=p^\text{sp}\sum_{q\in Q_i}\sum_{k\in K}{T_k\lambda_{q,i,k} \gamma_{q,i,k}}$, $\forall i \in \Psi_{}^{\text{tn}}$. As a result, the sub-problem (\ref{obj3})--(\ref{eqn:3_13}) is reduced to a simple optimal AC power flow problem. Based on the above assumption, there is a feasible solution $Y^{*}=\{l_{mn }, p_{i}^\text{ev}, p_{\text{un},i}^\text{ev}, s_{0 }, s_{m}^\text{pv}, s_{m}^\text{EV}, S_{mn},  v_{m }, v_m^\text{d}, \forall i \in \Psi_{}^{\text{tn}}, \forall m \in \Psi^{\text{dn}}, \forall \left(m,n\right) \in \Psi^{\text{db}} \} \in \text{relint}~\mathcal{D}$ subjects to:
\begin{align}
&|\underline{V_{m}}|^2 < v_{m } < |\overline{V_{m}}|^2, \qquad \forall m\in  \Psi^{\text{dn}}.
\end{align}

Furthermore, $\exists \Delta v>0$, subjects to:
\begin{align}
&|\underline{V_{m}}|^2 < v_{m } + \Delta v \leq  |\overline{V_{m}}|^2, \qquad \forall m\in  \Psi^{\text{dn}}.\label{eqn:deltav}
\end{align}

When $s_{m}^\text{pv} =0$, $\forall m \in \Psi^{\text{dn}}$, the active and reactive power injection at each node (except the root node $0$) are both negative. Therefore, the distribution system have nonzero unidirectional power flows so that we also have:
\begin{align}
&l_{mn}>0, \qquad\forall \{m,n\}\in \Psi^{\text{db}}.
\end{align}

There are only two non-affine constraints in each sub-problem, i.e., PV generation constraint (\ref{eqn:3_1}) and AC power flow constraint (\ref{eqn:3_7}). We discuss how we can construct a feasible solution based on $Y^{*}$ which strictly satisfies (\ref{eqn:3_1}) and (\ref{eqn:3_7}).

\subsubsection{PV Generation}\label{PVg}
In the first non-affine constraint (\ref{eqn:3_1}), the nameplate apparent power, i.e., $\overline{s_m^\text{pv}}, \forall m\in  \Psi^{\text{dn}}$, are nonnegative and given by the master problem. $\forall m\in  \Psi^{\text{dn}}$:
\begin{enumerate}[a)]
	\item If $\overline{s_m^\text{pv}}=0$, there is no PV generation at bus $m$ so that constraints (\ref{eqn:3_1})--(\ref{eqn:3_3}) can be omitted;
	\item Otherwise, $\overline{s_m^\text{pv}}>0$, for $s_{m}^\text{pv} =0$ in $Y^{*}$, it satisfies
	\begin{align}
	&\sqrt{|p_m^\text{pv}|^2+|q_m^\text{pv}|^2} = 0 < \overline{s_m^\text{pv}}.\label{eqn:PV1_3}
	\end{align}
\end{enumerate}

Therefore, $\forall m \in \Psi^{\text{dn}}$, if $\overline{s_m^\text{pv}}=0$, constraint (\ref{eqn:3_1}) can be omitted; otherwise, it is strictly satisfied for ${s_m^\text{pv}}=0$.

\subsubsection{AC Power Flow}

We slightly increase ${v}_m, \forall m \in \Psi^{\text{dn}}$, in $Y^{*}$ by the $\Delta v$ in constraint (\ref{eqn:deltav}) and adjust the corresponding nodal voltage deviations, i.e., ${v}_m^\text{d}, \forall m \in \Psi^{\text{dn}}$, to construct another solution
$Y^{**}=\{l_{mn }, p_{i}^\text{ev}, p_{\text{un},i}^\text{ev}, s_{0 }, s_{m}^\text{pv}, s_{m}^\text{EV}, S_{mn},  \hat{v}_{m}=v_{m }+\Delta v, \hat{v}_m^\text{d} = \max\{{v}_m^\text{d}, {v}_m^\text{d}+\Delta v\}, \forall i \in \Psi_{}^{\text{tn}}, \forall m \in \Psi^{\text{dn}}, \forall \left(m,n\right) \in \Psi^{\text{db}} \} \in \text{relint}~\mathcal{D}$. The other variables are equal to those in $Y^{*}$.
Then, we have 
\begin{align}
	& \forall m \in \Psi^{\text{dn}}, \forall \left(m,n\right) \in \Psi^{\text{db}}:\notag\\
	&\hat{v}_{m }-\hat{v}_{n } = {v}_{m} - {v}_{n }, \\
	&|\underline{V_{m}}|^2\leq \hat{v}_{m } = v_{m }+\Delta v\leq |\overline{V_{m}}|^2,\\
	&|S_{mn }|^2 \leq l_{mn }v_{m } <  l_{mn }(v_{m } + \Delta v) = l_{mn }\hat{v}_{m }, \label{acstrict}
\end{align}

As a result, the new solution $Y^{**}$ is still feasible and strictly satisfies the non-affine constraint (\ref{eqn:3_7}), i.e., (\ref{acstrict}). Besides, from Appendix \ref{PVg}, we also know that $Y^{**}$ strictly satisfies constraint (\ref{eqn:3_1}), when $\overline{s_m^\text{pv}}>0$.

To conclude, $Y^{**}\in \text{relint}~\mathcal{D}$ is a feasible solution for the sub-problem (\ref{obj3})--(\ref{eqn:3_13}), i.e., problem (\ref{sub_1})--(\ref{sub_2}), and it strictly satisfies all the non-affine constraints. Based on Proposition 1, we can conclude that strong duality holds between the sub-problem (\ref{sub_1})--(\ref{sub_2}) and its dual problem (\ref{subdual_1})--(\ref{subdual_3}). 

Moreover, because the total PV generation is constrained, the selling power of the system, i.e., $p_{0,\omega,t}^{-}$, is limited and the second term in (\ref{obj3}) is bounded below. The other terms in (\ref{obj3}) are all nonnegative. Thus, we can conclude that the sub-problem's objective (\ref{obj3}) is bounded below. Therefore, there exist a primal solution $Y^{**}$ and a dual solution $\{\mu^{*}, u^{*}\}$ that let the primal optimal objective $p^*$ equal to the dual optimal objective $q^*$.

\bibliographystyle{ieeetr}
\bibliography{ref}